\documentclass[amssymb, 11pt]{amsart}
\usepackage{latexsym}

\newlength{\standardunitlength}
\newtheorem{prop}{Proposition}[section]

\newtheorem{lemma}[prop]{Lemma}
\newtheorem{cor}[prop]{Corollary}
\newtheorem{theorem}[prop]{Theorem}

\begin{document}

\begin{center} \title [Martingales and character ratios]{\bf Martingales and character ratios} \end{center}

\author{Jason Fulman} \address{University of Pittsburgh\\ Pittsburgh,
PA 15260} \email{fulman@math.pitt.edu}



\date{February 19, 2004}

\begin{abstract} Some general connections between martingales and character ratios of finite groups are developed. As an application we sharpen the convergence rate in a central limit theorem for the character ratio of a random representation of the symmetric group on transpositions. A generalization of these results is given for Jack measure on partitions. We also give a probabilistic proof of a result of Burnside and Brauer on the decomposition of tensor products. \end{abstract}

\maketitle

\section{Introduction}

	The Plancherel measure of a finite group $G$ is a probability measure on the set of irreducible representations of $G$ which chooses a representation $\rho$ with probability
$\frac{dim(\rho)^2}{|G|}$, where $dim(\rho)$ denotes the dimension of $\rho$. For instance if $G$ is the symmetric group, the irreducible representations are parameterized by partitions $\lambda$ of $n$, and the
Plancherel measure chooses a partition $\lambda$ with probability $\frac{n!}{\prod_{x \in \lambda} h(x)^2}$ where
the product is over boxes in the partition and $h(x)$ is the
hooklength of a box. The hooklength of a box $x$ is defined as
1 + number of boxes in same row as x and to the right of x + number of boxes in same column of x and below x. For example
we have filled in each box in the partition of 7 below with
its hooklength \[ \begin{array}{c c c c} \framebox{6}&
	\framebox{4}& \framebox{2}& \framebox{1}\\ \framebox{3}&
	\framebox{1}&& \\ \framebox{1} &&& \end{array}, \] and the Plancherel measure would choose this partition with probability $\frac{7!}{(6*4*3*2)^2}$. Recently there has been
interest in the statistical properties of partitions chosen
from Plancherel measure. As it is beyond the scope of this paper to survey the topic, we refer the reader to the surveys
\cite{AD}, \cite{De} and the seminal papers \cite{J},
\cite{O1}, \cite{BOO} for a glimpse of the remarkable recent work on Plancherel measure.

	Let $\lambda$ be a partition of $n$ chosen from the Plancherel measure of the symmetric group $S_n$ and let $\chi^{\lambda}(12)$ be the irreducible character parameterized by $\lambda$ evaluated on the transposition $(12)$. The quantity $\frac{\chi^{\lambda}(12)}{dim(\lambda)}$ is called a character ratio and is crucial for analyzing the convergence rate of the random walk on the symmetric group generated by transpositions \cite{DS}. In fact Diaconis and Shahshahani prove that the eigenvalues for this random walk are the character ratios $\frac{\chi^{\lambda}(12)}{dim(\lambda)}$ each occurring with multiplicity $dim(\lambda)^2$. Character ratios on transpositions also play an essential role in work on the moduli space of curves \cite{EO}, \cite{OP}. 

	Given these motivations, it is natural to study the distribution of the character ratio $\frac{\chi^{\lambda}(12)}{dim(\lambda)}$. Kerov \cite{K1} outlined a proof of the following central limit theorem. 

\begin{theorem} \label{kerovs} (\cite{K1}) Let $\lambda$ be a partition of $n$ chosen from the Plancherel measure of the symmetric group $S_n$. Then the random variable $\frac{n-1}{\sqrt{2}} \frac{\chi^{\lambda}(12)}{dim(\lambda)}$ is asymptotically normal with mean 0 and variance 1. \end{theorem} 
	
	The details of the proof appeared in the paper \cite{IO},
	which gave a wonderful development of Kerov's work. Another
	proof, due to Hora \cite{Ho}, exploited the fact that the $k$th
	moment of a Plancherel distributed character ratio is the
	chance that the random walk generated by random transpositions
	is at the identity after k steps. Both of these proofs were
	essentially combinatorial in nature and used the method of
	moments (so information about all moments of the character
	ratio) to prove Theorem \ref{kerovs}.

	A more probabilistic approach to Kerov's central limit theorem appeared in \cite{F1}, which proved the following result. In its statement and throughout this paper, we use the notation that if $A$ is an event, $P(A)$ denotes the probability of $A$.
\begin{theorem} \label{my} (\cite{F1}) For $n \geq 2$ and all real $x_0$, \[ \left|P \left( \frac{n-1}{\sqrt{2}} \frac{\chi^{\lambda}(12)}{dim(\lambda)} \leq
x_0 \right) - \frac{1}{\sqrt{2 \pi}} \int_{-\infty}^{x_0}
e^{-\frac{x^2}{2}} dx \right| \leq \ 40.1 n^{-1/4}. \] \end{theorem} The proof technique involved random walk on the set of irreducible representations of the symmetric group and a technique known as Stein's method which is fundamentally different from the method of moments. Note that unlike Kerov's original result, Theorem \ref{my} included an error term (of order $n^{-1/4}$).

	In this paper we give a completely different approach to
	central limit theorems for character ratios. To begin we
	describe the general setting. Let $G$ be a finite group and
	$C$ a conjugacy class on which all irreducible characters of
	$G$ are real valued (this occurs for instance if
	$C=C^{-1}$). Let $\rho$ be a random representation of $G$,
	chosen from the Plancherel measure of $G$, and let
	$\chi^{\rho}$ denote the character of $\rho$. From the
	orthogonality relations of irreducible characters it follows that
	if $C$ is a nonidentity conjugacy class, then the random
	variable $\frac{|C|^{\frac{1}{2}} \chi^{\rho}(C)}{dim(\rho)}$
	has mean 0 and variance 1, and it is natural to investigate
	how close it is to a standard normal random variable with mean
	0 and variance 1.

	To see the connection with martingales, instead of
	investigating $\frac{|C|^{\frac{1}{2}}
	\chi^{\rho}(C)}{dim(\rho)}$, consider
	$\frac{|C|\chi^{\rho}(C)}{dim(\rho)}$. Given a subgroup
	chain $H_1=id \subseteq H_2 \subseteq \cdots \subseteq
	H_n=G$, together with representations $\rho(j) \in Irr(H_j)$,
	one can write \begin{eqnarray*} && \frac{|C|\chi^{\rho}(C)}{dim(\rho)}\\ & = &
	\sum_{j=2}^n \left( \frac{|C \cap H_j| \chi^{\rho(j)}(C \cap
	H_j)}{dim(\rho(j))} - \frac{|C \cap H_{j-1}|
	\chi^{\rho(j-1)}(C \cap H_{j-1})}{dim(\rho(j-1))} \right),\end{eqnarray*}
	where we set $\chi^{\rho(j)}(C \cap H_j)=0$ if $|C \cap
	H_j|=0$. Starting from $\rho(1)$ equal to the trivial
	representation of $H_1$, there is a way (see Section
	\ref{charmart} for details) to grow a series of random
	representations $\rho(2),\cdots,\rho(n)$ such that each
	$\rho(i)$ is chosen from the Plancherel measure of $H_i$. Thus
	$\frac{|C|\chi^{\rho}(C)}{dim(\rho)}$ is expressed as a sum of
	random variables, which one might hope to be small and nearly
	independent of each other, so that after renormalization a
	central limit theorem would hold.

Section \ref{charmart} proves that under certain conditions the sequence of partial sums in the above decomposition of $\frac{|C|\chi^{\rho}(C)}{dim(\rho)}$ is a martingale. It also shows that (due to the extra structure coming from representation theory), this martingale has special properties not shared by general martingales. There is a significant literature on central limit theorems for martingales, and in Section \ref{CLT} it is applied to prove the Theorem \ref{sharptrans} (actually we just prove Theorem \ref{sharptrans2} below, which generalizes Theorem \ref{sharptrans}).

\begin{theorem} \label{sharptrans} For any $s<\frac{1}{2}$, there is a constant $A_s$ so that for all $n$ and real $x_0$, \[ \left|P \left( \sqrt{{n-1 \choose 2}} \frac{\chi^{\lambda}(12)}{dim(\lambda)} \leq x_0 \right) - \frac{1}{\sqrt{2 \pi}} \int_{-\infty}^{x_0}
e^{-\frac{x^2}{2}} dx \right| \leq \ A_s n^{-s}. \] \end{theorem} This sharpens the convergence rate in Theorem \ref{my} and very nearly proves the conjecture in \cite{F1} that in Theorem \ref{sharptrans} one could take $s=\frac{1}{2}$. An essential ingredient in the proof of Theorem \ref{sharptrans} is an understanding of the conditional moments of the summands in the decomposition of $\frac{|C|\chi^{\rho}(C)}{dim(\rho)}$. For general conjugacy classes this is a difficult problem. A fortunate and remarkable fact is that work of Kerov \cite{K3} in the combinatorics literature implies that asymptotically these conditional moments are simply the moments of the semicircle distribution. Section \ref{trans} explains this, and discusses some other phenomena which occur only for the class of transpositions. 

	It should be pointed out that both Kerov \cite{K1} and Hora \cite{Ho} use the method of moments to prove central limit theorems (with no error term) for character ratios evaluated on cycles of various lengths-even obtaining a multidimensional central limit theorem showing these to be asymptotically independent. Martingales do blend well with multidimensional central limit theorems, but it is unclear whether the methods of this paper can be pushed through. Hora \cite{Ho} also shows that there are many conjugacy classes of the symmetric group where central limit theorems do not hold.

	This paper also considers central limit theorems for
Jack$_{\alpha}$ measure, a generalization of Plancherel measure of the
symmetric group on partitions of size $n$. Here $\alpha>0$, and this
measure chooses a partition $\lambda$ of size $n$ with probability
\[\frac{\alpha^n n!}{\prod_{s \in \lambda} (\alpha a(s) + l(s) +1)
(\alpha a(s) + l(s) + \alpha)}, \] where the product is over all boxes
in the partition. Here $a(s)$ denotes the number of boxes in the same
row of $s$ and to the right of $s$ (the ``arm'' of s) and $l(s)$
denotes the number of boxes in the same column of $s$ and below $s$
(the ``leg'' of s). For example the partition of 5 below \[
\begin{array}{c c c} \framebox{\ }& \framebox{\ }& \framebox{\ } \\
\framebox{\ }& \framebox{\ }& \end{array}, \] would have Jack$_{\alpha}$
measure \[\frac{60 \alpha^2}{(2 \alpha+2)(3 \alpha+1) (\alpha+2)(2
\alpha+1)(\alpha+1)}.\] Note that when $\alpha=1$, Jack measure
reduces to Plancherel measure of the symmetric group. The papers
\cite{O2}, \cite{BO1} emphasize that for $\alpha$ fixed the study of Jack$_{\alpha}$ measure is an important open problem, about which relatively little is known for general values of $\alpha$ (the three
values $\alpha=1/2,1,2$ are tractable). It is a discrete analog of ensembles from random matrix theory and like Jack polynomials
\cite{GHJ}, should also be relevant to the moduli space of curves.

	Recently the author proved a central limit related to Jack measure, for general values of $\alpha$. We remark that when $\alpha=2$ the statistic $W_{\alpha}$ is a spherical function for the Gelfand pair $(S_{2n},H_{2n})$, where $S_{2n}$ is a symmetric group and $H_{2n}$ is the hyperoctahedral group of size $2^n n!$.

\begin{theorem} \label{my2} (\cite{F2}) Suppose that $\alpha \geq 1$. Let $W_{\alpha}(\lambda) = \frac{\alpha n(\lambda')-n(\lambda)}{\sqrt{\alpha {n \choose 2}}}$. For $n \geq 2$ and all real $x_0$, \[ \left|P_{\alpha} \left( W_{\alpha} \leq
x_0 \right) - \frac{1}{\sqrt{2 \pi}} \int_{-\infty}^{x_0}
e^{-\frac{x^2}{2}} dx \right| \leq A_{\alpha} n^{-1/4} \] where $A_{\alpha}$ depends on $\alpha$ but not on $n$.
\end{theorem}
 
	Section \ref{jack} develops connections between Jack
symmetric functions and martingales, showing how statistics such as $W_{\alpha}$ can be written as a sum of martingale
differences. Although the results are analogous to those of
Section \ref{charmart} (and generalize them for the case of the symmetric group), the proofs are completely different as
representation theory of finite groups can no longer be
used. Section \ref{trans2} then focuses on the summands of the specific statistic $W_{\alpha}$. Again one finds that certain simplifications occur which do not hold in general for the martingales in Section \ref{jack}. In Section \ref{CLT} the results of Sections \ref{jack} and \ref{trans2} are used to prove the following result. 

\begin{theorem}
 \label{sharptrans2} Suppose that $\alpha \geq 1$. Let
$W_{\alpha}(\lambda) = \frac{\alpha
n(\lambda')-n(\lambda)}{\sqrt{\alpha {n \choose 2}}}$. Then for any
$s<\frac{1}{2}$, there is a constant $A_{s,\alpha}$ (depending only on
$s$ and $\alpha$) such that for all $n \geq 1$ and real $x_0$, \[
\left|P_{\alpha} \left( W_{\alpha} \leq x_0 \right) - \frac{1}{\sqrt{2
\pi}} \int_{-\infty}^{x_0} e^{-\frac{x^2}{2}} dx \right | \leq
A_{s,\alpha} n^{-s} .\]
\end{theorem} 

The assumption that $\alpha \geq 1$ in Theorems \ref{my2} and \ref{sharptrans2} is merely for convenience. Indeed, from the definition of Jack measure it is clear that the Jack$_{\alpha}$ probability of $\lambda$ is equal to the Jack$_{\frac{1}{\alpha}}$ probability of the transpose of $\lambda$. From this one concludes that the Jack$_{\alpha}$ probability that $W_{\alpha}=w$ is equal to the Jack$_{\frac{1}{\alpha}}$ probability that $W_{\frac{1}{\alpha}}=-w$, so that a central limit theorem holds for $\alpha$ if and only if it holds for $\frac{1}{\alpha}$.

{\bf Conjecture 1:} Theorem \ref{sharptrans2} also holds with $s=\frac{1}{2}$.

	Conjecture 1 was stated in \cite{F1} for Plancherel measure. Aside from Theorem \ref{sharptrans2} itself, further evidence for Conjecture 1 is Proposition 6.5 of \cite{F2}. 

	Given that the convergence rates in Theorems \ref{my} and
	\ref{my2} can be improved by martingale theory, it is natural
	to ask if the ``exchangeable pairs'' Stein's method techniques
	of those papers are obsolete. The answer is a firm no. First,
	exchangeable pairs will be crucial for a Stein's method
	approach to other statistics of Plancherel measure and Jack
	measure. Second, the construction of exchangeable pairs
	involved random walk on the set of irreducible representations
	of a finite group and was shown in \cite{F1},\cite{F2} to have
	some implications for the decomposition of tensor products. In
	Section \ref{probproof} of this paper, it is shown how these
	random walks give a probabilistic proof of a theorem of
	Burnside and Brauer (page 49 of \cite{Is}) which states that
	if a faithful character $\eta$ of a finite group assumes $m$
	distinct values, then every irreducible character of $G$
	occurs as a constituent of some $\eta^j$, where $0 \leq
	j<m$. Although the proof is less general than their result (it
	only works for real valued characters), the conceptual
	framework is more appealing.

	The organization of this paper is as follows. Section
	\ref{charmart} develops connections between martingales and
	character ratios of finite groups. Section \ref{trans} then
	focuses on the case of transpositions, where certain
	simplifications occur. Section \ref{jack} develops general
	connections between Jack symmetric functions and martingales
	and is highly combinatorial in nature. Section \ref{trans2}
	then focuses on combinatorics related to the statistic
	$W_{\alpha}$, where again there are attractive
	simplifications. (The reader will observe that parts of
	Section \ref{trans2} generalize parts of Section \ref{trans},
	but there are real differences and since both sections are
	brief and Section \ref{trans} might be unreadable by group
	theorists who work on character ratios but not on Jack
	polynomials, they have been separated). Section \ref{CLT} uses
	central limit theory for martingales to prove Theorem
	\ref{sharptrans2} (and so also Theorem
	\ref{sharptrans}). Finally, Section \ref{probproof} gives a
	probabilistic proof of a result of Burnside and Brauer on the
	decomposition of tensor products, thereby showing that the
	machinery underlying a Stein's method approach to Theorems
	\ref{my} and \ref{my2} is of ongoing interest.

\section{Character ratios and martingales} \label{charmart}

	Let $G$ be a finite group and $C$ a conjugacy class of $G$ on which all irreducible characters of $G$ are real (this holds for instance if $C=C^{-1}$). Suppose there is a chain of subgroups $H_1=id \subseteq H_2 \subseteq \cdots \subseteq H_n=G$. Given a representation $\rho(j)$ of $H_j$, one obtains a random representation $\rho(j+1)$ of $H_{j+1}$ by choosing $\rho(j+1)$ with probability \[ \frac{|H_j| dim(\rho(j+1)) \kappa(\rho(j),\rho(j+1))}{|H_{j+1}| dim(\rho(j))}. \] Here $\kappa(\rho(j),\rho(j+1))$ denotes the multiplicity of $\rho(j+1)$ in the induction of $\rho(j)$ from $H_{j}$ to $H_{j+1}$, or by Frobenius reciprocity the multiplicity of $\rho(j)$ in the restriction of $\rho(j+1)$ from $H_{j+1}$ to $H_j$. Thus starting from $\rho(1)$ (the trivial representation of a one element group) one obtains a random sequence $(\rho(1),\cdots,\rho(n))$ of representations. It is proved in \cite{F2} that each $\rho(j)$ is distributed according to the Plancherel measure of $H_j$ (this was known earlier for the symmetric group \cite{K2}).

	A sequence $(Y_1,\cdots,Y_n)$ of real valued random variables
	is called a martingale with respect to the sequence of random
	variables $(X_1,\cdots,X_n)$ if two properties hold. First, for
	each $n \geq 1$ there is a function $f_n$ such that
	$Y_n=f_n(X_1,\cdots,X_n)$. Second, the martingale identity \[
	E(Y_j|X_1,\cdots,X_{j-1})= Y_{j-1} \] must hold for all $j \geq
	1$. Throughout this section we also define $Y_0$ to be
	identically equal to 0.

\begin{theorem} \label{ismart} Let $G$ be a finite group and $C=C^{-1}$ a
 self-inverse conjugacy class of $G$. Let $\rho(1),\cdots,\rho(n)$ be
 the sequence of random representations associated to a subgroup chain
 $H_1=id \subseteq \cdots \subseteq H_n=G$. Suppose that $C_j=C \cap
 H_j$ is either empty or a single conjugacy class of $H_j$. Let
 $Y_j=\frac{|C_j| \chi^{\rho(j)}(C_j)}{dim(\rho(j))}$, where this is
 understood to be 0 if $|C_j|=0$. Then the sequence $(Y_1,\cdots,Y_n)$
 is a martingale. \end{theorem}

\begin{proof} Let $Irr(H)$ denote the set of irreducible representations of a finite group $H$. Observe that the expected value of $Y_{j+1}$ given the earlier $\rho$'s is \begin{eqnarray*}
& & \sum_{\tau \in Irr(H_{j+1})} \frac{|H_j| dim(\tau)
\kappa(\rho(j),\tau)}{|H_{j+1}| dim(\rho(j))}\frac{|C_{j+1}|
\chi^{\tau}(C_{j+1})}{dim(\tau)}\\ & = & \frac{|C_{j+1}|
|H_j|}{dim(\rho(j)) |H_{j+1}|} \sum_{\tau \in Irr(H_{j+1})}
\kappa(\rho(j),\tau) \chi^{\tau}(C_{j+1})\\ & = & \frac{|C_{j+1}|
|H_j|}{dim(\rho(j)) |H_{j+1}|}
Ind^{H_{j+1}}_{H_j}(\chi^{\rho(j)})[C_{j+1}]. \end{eqnarray*} Here
$Ind$ stands for induction of characters. If $C_{j+1}$ is empty, then
$Y_{j+1}=Y_j=0$. Otherwise since $C_j$ is either empty or a single conjugacy class of $H_j$, it follows from a general formula
for induced characters (page 30 of \cite{Se}) that \[
Ind^{H_{j+1}}_{H_j}(\chi^{\rho(j)})[C_{j+1}] =
\frac{|H_{j+1}||C_j|}{|H_j||C_{j+1}|} \chi^{\rho(j)}(C_j).\]
Substituting this in implies that the expected value of $Y_{j+1}$
given the earlier $\rho$'s is equal to $Y_j$.\end{proof}

{\bf Remarks:}
\begin{enumerate}
\item The conditions of Theorem \ref{ismart} may seem restrictive but
do apply in many cases of interest, most notably the symmetric groups
and for many classes in finite classical groups.
\item Note also that from page 52 of \cite{Se} it follows that the random variables $Y_1,\cdots,Y_n$ are algebraic integers. Since the characters of the symmetric groups are integer valued, in this case $Y_1,\cdots,Y_n$ must be rational algebraic integers, and hence integers.
\end{enumerate}

	Proposition \ref{condit} computes certain conditional probabilities related to the martingale sequence $(Y_1,\cdots,Y_n)$, which is not possible for general martingales. It will not be needed in the proof of central limit theorems. We use the notation that $E^W(Z)$ denotes the expected value of $Z$ given $W$.

\begin{prop} \label{condit} Assume the hypotheses of Theorem \ref{ismart}. Then for all $1 \leq j \leq n$, \[ E^{Y_n}(Y_j) = \frac{|C_j|}{|C_n|} Y_n.\] \end{prop}

\begin{proof} In fact we show that \[ E^{\rho(n)}(Y_j) = \frac{|C_j|}{|C_n|} \frac{|C_n| \chi^{\rho(n)}(C_n)}{dim(\rho(n))} \] which is stronger. 

	For $j=n$ this is clear, so consider $j=n-1$. For ease of notation we let $\lambda$ denote $\rho(n)$. Then
\[ E^{\rho(n)}(Y_{n-1}) =  \sum_{\tau \in Irr(H_{n-1})}
\frac{P(\rho(n-1)=\tau,\rho(n)=\lambda)}{P(\rho(n)=\lambda)}
\frac{|C_{n-1}| \chi^{\tau}(C_{n-1})}{dim(\tau)}.\] Now since $\rho(n)$ has the distribution of Plancherel measure of $H_n$, $P(\rho(n)=\lambda)=\frac{dim(\lambda)^2}{|H_n|}$. Similarly
\begin{eqnarray*} P(\rho(n-1)=\tau,\rho(n)=\lambda) &=&
\frac{dim(\tau)^2}{|H_{n-1}|} \frac{\kappa(\tau,\lambda) dim(\lambda)|H_{n-1}|}{dim(\tau) |H_n|}\\ & = & \frac{dim(\tau) dim(\lambda) \kappa(\tau,\lambda)}{|H_n|}. \end{eqnarray*} Substituting this into the expression for $ E^{\rho(n)}(Y_{n-1})$ shows that as desired, the latter is equal to \[ \frac{|C_{n-1}|}{dim(\lambda)} \sum_{\tau \in Irr(H_{n-1})} \kappa(\tau,\lambda)
\chi^{\tau}(C_{n-1}) = \frac{|C_{n-1}|
\chi^{\lambda}(C_n)}{dim(\lambda)}.\]

	Next we indicate how to treat the case $j=n-2$, from which the general argument is completely clear. By considerations similar to those in the previous paragraph, $E^{\rho(n)}(Y_{n-2})$ is equal to \begin{eqnarray*} & & \sum_{\tau \in Irr(H_{n-2})}
\frac{P(\rho(n-2)=\tau,\rho(n)=\lambda)}{P(\rho(n)=\lambda)}
\frac{|C_{n-2}| \chi^{\tau}(C_{n-2})}{dim(\tau)}\\
 & = & \sum_{\tau \in
Irr(H_{n-2})} \sum_{\eta \in Irr(H_{n-1})}
\frac{P(\rho(n-2)=\tau,\rho(n-1)=\eta,\rho(n)=\lambda)}{ P(\rho(n)=\lambda) }\\
& & \cdot
\frac{|C_{n-2}| \chi^{\tau}(C_{n-2})}{dim(\tau)}\\
 & = & \frac{|C_{n-2}|}{dim(\lambda)} \sum_{\tau \in Irr(H_{n-2})} \sum_{\eta
\in Irr(H_{n-1})} \kappa(\tau,\eta) \kappa(\eta,\lambda)
\chi^{\tau}(C_{n-2})\\
 & = & 
\frac{|C_{n-2}|\chi^{\lambda}(C_{n})}{dim(\lambda)}. \end{eqnarray*}
Note that the last equality used the fact that restriction of
characters is transitive (i.e. to compute the restriction of $\chi$ to elements in $H_{n-2}$, one can first restrict to $H_{n-1}$ and then to $H_{n-2}$). \end{proof}

	Proposition \ref{exactfor} derives an exact expression (which will be needed) for $E(Y_j-Y_{j-1})^2$.

\begin{prop} \label{exactfor} Assume the hypotheses of Theorem \ref{ismart}. Then for $j \geq 2$, $E(Y_j-Y_{j-1})^2=|C_j|-|C_{j-1}|$. \end{prop} 

\begin{proof} Clearly $E(Y_j-Y_{j-1})^2=E(Y_j^2)+E(Y_{j-1}^2)-2E(Y_jY_{j-1})$. But \[ E(Y_j Y_{j-1}) = E(E^{Y_{j-1}}(Y_jY_{j-1})) = E(Y_{j-1}E^{Y_{j-1}}(Y_j)) = E(Y_{j-1}^2).\] Thus $E(Y_j-Y_{j-1})^2$ is equal to $E(Y_j^2)-E(Y_{j-1}^2)$. Now observe that \[ E(Y_j^2) = \sum_{\tau \in Irr(H_j)} \frac{dim(\tau)^2}{|H_j|} \frac{|C_j|^2 \chi^{\tau}(C_j)^2}{dim(\tau)^2}.\] After canceling the factors of $dim(\tau)^2$, it follows from the orthogonality relations of the characters of $H_j$ that $E(Y_j^2)=|C_j|$ for all $j$. This implies the proposition. \end{proof}

\section{Transpositions} \label{trans}

	This section illustrates and sharpens the results in Section \ref{charmart} for the case of character ratios on
transpositions. Along the way we make some links with the
combinatorics literature. Moreover certain properties will be established for the class of transpositions which fail for other conjugacy classes; for instance it will be proved that not only the sequence $(Y_1,\cdots,Y_n)$, but also the sequence \[ (Y_1^2,Y_2^2-{2 \choose 2}, \cdots, Y_n^2-{n \choose 2})\] is a martingale.

	Throughout this section $G=S_n$, $H_j=S_j$ for $1 \leq j \leq n$, and $C$ is the conjugacy class of transpositions. Recall that the irreducible representations of the symmetric group are parameterized by partitions $\lambda$ of $n$ (see \cite{Sa} for a friendly introduction to representation theory of the symmetric group). Frobenius \cite{Fr} found the following explicit formula for the character ratio of the symmetric group on transpositions: \[ \frac{\chi^{\lambda}(12)}{dim(\lambda)} = \frac{1}{{n \choose 2}} \sum_i \left({\lambda_i \choose 2} - {\lambda_i' \choose 2} \right)\] where $\lambda_i$ is the length of row $i$ of $\lambda$ and $\lambda_i'$ is the length of column $i$ of $\lambda$. 

From the branching rules of the symmetric group, the multiplicity $\kappa(\mu,\lambda)$ is non-negative (and in fact is 1) if and only $\lambda$ is obtained from $\mu$ by adding a box to the diagram of $\mu$. From Frobenius' formula it is follows that if $\mu$ is such a partition, then \[ \frac{|C_n| \chi^{\lambda}(12)}{dim(\lambda)}- \frac{|C_{n-1}| \chi^{\mu}(12)}{dim(\mu)} =c(x) \] where $x$ is the box added to $\mu$ to obtain $\lambda$ and the content of a box $c(x)$ is defined as column number of box - row number of box.

Given the above discussion, one sees that Theorem \ref{ismart} implies that for any partition $\lambda=\rho(j)$ of size $j$, the expected value of the content of the box added to $\lambda$ in moving from $\rho(j)$ to $\rho(j+1)$ is 0. Rather remarkably this combinatorial fact was known to Kerov \cite{K2},\cite{K3}, who had also proved that for {\it any} $\lambda$ of size $j$, the variance of the content of the box added to $\lambda$ is equal to $j$. Note that this is stronger than Proposition \ref{exactfor}, which only establishes this for a random Plancherel distributed $\lambda$. For conjugacy classes other than transpositions, it need not be the case that $E^{\rho(j)}(Y_{j+1}-Y_{j})^2$ depends only on $\rho(j)$. Indeed the reader can check that this fails for the class of 3-cycles and $j=3$. 

For the purpose of proving Theorem \ref{sharptrans}, it will be necessary to have upper bounds on all moments of the variable $E(Y_{j+1}-Y_j)$. In very recent work Lassalle found formulas for all $E^{\rho(j)}(Y_{j+1}-Y_{j})^r$, from which one can obtain $E(Y_{j+1}-Y_j)^r$ by averaging with respect to Plancherel measure. The formulas are quite complicated but for instance for $r=1,2,3,4$ one has the following special case of Theorem 8.1 of \cite{La}.
\begin{theorem} \label{Lmoment} (\cite{La}) Let $d_k(\lambda)$ denote $\sum_{x \in \lambda} c(x)^k$. Let $s_r(\lambda)$ be the $r$th moment of the content of the box added to $\rho(j)=\lambda$ when transitioning to $\rho(j+1)$. Then \[ \begin{array}{ll}
s_1(\lambda)=0 & \\
s_2(\lambda)=|\lambda| & \\
s_3(\lambda)= 2 d_1(\lambda) & \\
s_4(\lambda)= 3 d_2(\lambda) + {|\lambda|+1 \choose 2} &
\end{array}\] \end{theorem} 

	In Section \ref{trans2} it is shown (in the Jack setting which is more general) how Lassalle's formulas lead to the necessary upper bounds on $E(Y_{j+1}-Y_j)^r$. But in the case of Plancherel measure the exact asymptotics of $E(Y_{j+1}-Y_j)^r$ is known. Indeed, one obtains simply the moments of the semi-circle distribution. We only state the result for even moments, since the odd moments all vanish by symmetry of Plancherel measure under transposing the partition. Theorem \ref{ex} is due to Kerov \cite{K3}. It also follows in a sentence from work of Biane \cite{Bia} provided that one is familiar with Murphy's work \cite{Mu} relating the spectrum of the element $(1,2)+\cdots+(1,n)$ to the content of box number $n$ in standard Young tableaux.

\begin{theorem} \label{ex} (\cite{K3}) Let $\lambda$ be chosen from Plancherel measure on partitions of size $n$. Then \[ lim_{n \rightarrow \infty} \frac{E(s_{2r}(\lambda))}{n^r} = \frac{{2r \choose r}}{r+1}.\] \end{theorem}	

	The final results of this section are not needed for proving central limit theorems, but are of interest. Corollary \ref{record1} shows that for the class of transpositions, Proposition \ref{condit} has a curious combinatorial consequence. We remind the reader that a standard Young tableau of shape $\lambda$ is a filling of the boxes of $\lambda$ with the numbers $1,\cdots,|\lambda|$ each occurring once such that the numbers increase going across rows (from left to right) and down columns.

\begin{cor} \label{record1} Let $T$ be a standard Young tableau of shape $\lambda$ chosen uniformly at random. Then the expected value of the content of box $j$ is equal to $\frac{j-1}{{n \choose 2}} \sum_{x \in \lambda} c(x)$. \end{cor}

\begin{proof} From the discussion at the beginning of this section, $Y_j-Y_{j-1}$ is the content of box $j$. By the proof of Proposition \ref{condit}, the expected value of
 $Y_j-Y_{j-1}$ given that $\rho(n)=\lambda$ is equal to \[ \frac{{j \choose 2}-{j-1 \choose 2}}{{n
 \choose 2}} \sum_{x \in \lambda} c(x)= \frac{j-1}{{n
 \choose 2}} \sum_{x \in \lambda} c(x).\] Now observe that given that
 $\rho(n)=\lambda$, all sequences $(\rho(1),\cdots,\rho(n)=\lambda)$
 occur with probability \[ \frac{1}{dim(\lambda)^2/n!}
 \prod_{j=1}^{n-1} \frac{j! dim(\rho(j+1))}{(j+1)! dim(\rho(j))} =
 \frac{1}{dim(\lambda)}.\] Thus such sequences correspond to standard
 Young tableaux $T$ under the uniform distribution. \end{proof}

	Finally we prove that for the conjugacy class of
	transpositions the sequence $(Y_1^2,Y_2^2-{2 \choose 2},
	\cdots, Y_n^2-{n \choose 2})$ is a martingale.

\begin{prop} \label{squaremart} Let $C$ be the conjugacy class of transpositions in $S_n$, and let $(Y_1,\cdots,Y_n)$ be the martingale of Theorem \ref{ismart}. Then $(Y_1^2,Y_2^2-{2 \choose 2}, \cdots, Y_n^2-{n \choose 2})$ is also a martingale.
\end{prop}

\begin{proof} Observe that for all $\mu$ of size $j$, \begin{eqnarray*} & &  - {j+1 \choose 2} + \sum_{|\tau|=j+1} P(\rho(j+1)=\tau|\rho(j)=\mu) \left (\sum_{x \in \tau} c(x) \right)^2 \\
& = & -{j+1 \choose 2} + \sum_{|\tau|=j+1} P(\rho(j+1)=\tau|\rho(j)=\mu)\\
& & \cdot \left((\sum_{x \in \mu} c(x))^2 + c(y)^2 + 2 c(y) ( \sum_{x \in \mu} c(x)) \right), \end{eqnarray*} where $y$ is the box of $\tau$ not in $\mu$. Since the expected value of $c(y)$ given $\mu$ is 0 and the expected value of $c(y)^2$ given $\mu$ is $j$, it follows that the last expression is equal to $Y_j^2  -{j \choose 2}$, as desired. \end{proof}

\section{Jack polynomials and martingales} \label{jack}

	The purpose of this section is to extend to Jack measure the connections with martingales in Section \ref{charmart}. The arguments are (by necessity) combinatorial as opposed to the algebraic arguments used in Section \ref{charmart}, and we assume that the reader is familiar with symmetric functions as in \cite{Mac}. 

	As in the introduction, given a box $s$ in the diagram of
	$\lambda$, let $a(s)$ and $l(s)$ denote the arm and leg of $s$
	respectively. One defines quantities \[ c_{\lambda}(\alpha) =
	\prod_{s \in \lambda} (\alpha a(s) + l(s) +1) \] \[
	c_{\lambda}'(\alpha) = \prod_{s \in \lambda} (\alpha a(s) +
	l(s) + \alpha) .\] Let $m_i(\lambda)$ be the number of
	parts (i.e. rows) of $\lambda$ of size $i$ and let
	$l(\lambda)$ denote the total number of parts of
	$\lambda$. The symbol $z_{\lambda}$ denotes $\prod_{i \geq 1}
	i^{m_i(\lambda)} m_i(\lambda)!$, the size of the centralizer
	of a permutation of cycle type $\lambda$ in the symmetric
	group.

	Let $\theta^{\lambda}_{\mu}(\alpha)$ be the coefficient of the power sum symmetric function $p_{\mu}$ in $J_{\lambda}^{(\alpha)}$. These will be the analogs of $\frac{|C| \chi^{\lambda}(C)}{dim(\lambda)}$ studied in Section \ref{charmart}. In fact when $\alpha=1$ they specialize to $\frac{|C| \chi^{\lambda}(C)}{dim(\lambda)}$, where $C$ is a conjugacy class of the symmetric group of type $\mu$. When $\alpha=2$ the $\theta^{\lambda}_{\mu}(\alpha)$ are spherical functions for the Gelfand pair $(S_{2n},H_{2n})$ where $S_{2n}$ is a symmetric group and $H_{2n}$ is the hyperoctahedral group of size $2^n n!$.

	Next we consider the ring of symmetric functions, with inner product defined by the orthogonality condition $<p_{\nu},p_{\mu}>_{\alpha}=\delta_{\nu,\mu} z_{\mu} \alpha^{l(\mu)}$. For a symmetric function $f$, its adjoint $f^{\perp}$ is defined by the condition $<fg,h>_{\alpha}=<g, f^{\perp}h>_{\alpha}$ for all $g,h$ in the ring of symmetric functions. It is straightforward to check using the basis of power sum symmetric functions that $p_1^{\perp}= \alpha \frac{\partial}{\partial p_1}$ (for the case $\alpha=1$ see page 76 of \cite{Mac}).

	Let \[ \psi_{\lambda/\tau}' = \prod_{s \in
C_{\lambda/\tau}-R_{\lambda/\tau}} \frac{(\alpha a_{\lambda}(s) +
l_{\lambda}(s)+1)}{(\alpha a_{\lambda}(s) + l_{\lambda}(s)+\alpha)}
\frac{(\alpha a_{\tau}(s) + l_{\tau}(s)+\alpha)}{(\alpha a_{\tau}(s) +
l_{\tau}(s)+1)} \] where $C_{\lambda/\tau}$ is the union of columns of $\lambda$ that intersect $\lambda-\tau$ and $R_{\lambda/\tau}$ is the union of rows of $\lambda$ that intersect $\lambda-\tau$.

	The following lemma will be useful.

\begin{lemma} \label{useful}
\begin{enumerate}
\item \cite{Mac} \[ \sum_{|\rho|=n} \frac{\theta^{\rho}_{\mu}(\alpha) \theta^{\rho}_{\eta}(\alpha)}{c_{\rho}(\alpha) c'_{\rho}(\alpha)} = \delta_{\mu,\eta} \frac{1}{z_{\mu} \alpha^{l(\mu)}}.\]
\item \cite{Mac} \[ p_1 J_{\lambda}^{(\alpha)} = \sum_{|\Lambda|=|\lambda|+1} \frac{c_{\lambda}(\alpha)}{c_{\Lambda}(\alpha)} \psi_{\Lambda/\lambda}'(\alpha) J_{\Lambda}^{(\alpha)}. \]
\item \cite{F2} $p_1^{\perp} J_{\lambda}^{(\alpha)} = \sum_{|\tau|=|\lambda|-1} \frac{c_{\lambda}'(\alpha) \psi'_{\lambda/\tau}(\alpha)}{c_{\tau}'(\alpha)} J_{\tau}^{(\alpha)}$.
\end{enumerate}
\end{lemma}

	Using Lemma \ref{useful} we establish some relations concerning the $\theta^{\lambda}_{\mu}(\alpha)$'s. We use the notation that $\mu+1^b$ denotes the partition of $|\mu|+b$ given by adding b parts of size 1 to $\mu$, and that $\mu-1^b$ denotes the partition of $|\mu|-b$ given by removing b parts of size 1 from $\mu$ (if this is possible). If $\mu$ has fewer than $b$ parts of size 1, we define $\theta^{\lambda}_{\mu-1^b}(\alpha)$ to be 0 for all $\lambda$. We let $m_1(\mu)$ denote the number of parts of $\mu$ of size 1.

\begin{lemma} \label{crux} 
\begin{enumerate} 
\item Let $|\lambda|=n$ and $|\mu|=n+1$. Then \[ \theta^{\lambda}_{\mu-1}(\alpha) = \sum_{|\Lambda|=n+1} \frac{c_{\lambda}(\alpha)}{c_{\Lambda}(\alpha)} \psi_{\Lambda/\lambda}'(\alpha) \theta^{\Lambda}_{\mu}(\alpha).\]
\item Let $|\lambda|=|\mu|=n$. Then \[ \alpha m_1(\mu)  \theta^{\lambda}_{\mu}(\alpha) = \sum_{|\tau|=n-1} \frac{c_{\lambda}'(\alpha)}{c_{\tau}'(\alpha)} \psi'_{\lambda/\tau}(\alpha) \theta^{\tau}_{\mu-1}(\alpha).\]
\end{enumerate}
\end{lemma}

\begin{proof} For the first assertion, consider the inner product $<p_1 J_{\lambda}^{(\alpha)}, p_{\mu}>$. On one hand, by part 2 of Lemma \ref{useful} it is equal to \begin{eqnarray*} & &
<\sum_{|\Lambda|=n+1} \frac{c_{\lambda}(\alpha)}{c_{\Lambda}(\alpha)}
\psi_{\Lambda/\lambda}'(\alpha) J_{\Lambda}^{(\alpha)}, p_{\mu}>\\ & =
& \sum_{|\Lambda|=n+1} \frac{c_{\lambda}(\alpha)}{c_{\Lambda}(\alpha)}
\psi_{\Lambda/\lambda}'(\alpha) \theta^{\Lambda}_{\mu}(\alpha) z_{\mu}
\alpha^{l(\mu)}. \end{eqnarray*} On the other hand it is equal to \[ <J_{\lambda}^{(\alpha)}, p_1^{\perp}
p_{\mu}> = <J_{\lambda}^{(\alpha)}, \alpha \frac{\partial}{\partial
p_1} p_{\mu}>.\] If $m_1(\mu)=0$ this is 0 so the assertion is
true. Otherwise it is \[ <J_{\lambda}^{(\alpha)}, \alpha(m_1(\mu))
p_{\mu-1}> = \alpha(m_1(\mu)) z_{\mu-1} \alpha^{l(\mu)-1}
\theta^{\lambda}_{\mu-1}(\alpha).\] Comparing these two expressions
and using the fact that $z_{\mu}=(m_1(\mu))z_{\mu-1}$ proves the first
assertion.

	For the second assertion, consider the inner product
	$<p_1^{\perp} J_{\lambda}^{(\alpha)},p_{\mu-1}>$. On one
	hand, by part 3 of Lemma \ref{useful} it is equal to \begin{eqnarray*} & &
	<\sum_{|\tau|=n-1}
	\frac{c_{\lambda}'(\alpha)}{c_{\tau}'(\alpha)}
	\psi_{\lambda/\tau}'(\alpha) J_{\tau}^{(\alpha)},p_{\mu-1}>\\ & =&
	\sum_{|\tau|=n-1}
	\frac{c_{\lambda}'(\alpha)}{c_{\tau}'(\alpha)}
	\psi_{\lambda/\tau}'(\alpha) \theta^{\tau}_{\mu-1}(\alpha)
	z_{\mu-1} \alpha^{l(\mu)-1}. \end{eqnarray*} On the other hand, it is equal
	to \[
	<J_{\lambda}^{(\alpha)},p_{\mu}>=\theta^{\lambda}_{\mu}(\alpha)
	z_{\mu} \alpha^{l(\mu)}.\] Comparing these two expressions and
	using the fact that $z_{\mu}=m_1(\mu) z_{\mu-1}$ proves the
	second assertion. \end{proof}

	In order to have analogs of the results of Section \ref{charmart}, it is necessary to have an analog of the growth process for representations of the symmetric group. Fortunately, such a process has been developed by Kerov \cite{K4}. His process is best understood in terms of harmonic functions on Bratelli diagrams, but rather than going into this (see \cite{K4} or Sections 2.1 and 4 of \cite{F2} for details) we simply describe in the next paragraph the process and state the properties we need.

It is convenient to define $dim_{\alpha}(\lambda)=\frac{n! \alpha^n}{c_{\lambda}'(\alpha)}$, which in the case $\alpha=1$ reduces to the dimension of the irreducible representation of the symmetric group parameterized by $\lambda$. A result of Stanley \cite{Sta} is that $dim_{\alpha}(\lambda)= \sum_{|\tau|=n-1} \psi'_{\lambda/\tau}(\alpha) dim_{\alpha}(\tau)$. In Kerov's growth process the chance of transitioning from a partition $\lambda$ of size $n$ to a partition $\Lambda$ of size $n+1$ is equal to \[ \frac{\psi_{\Lambda/\lambda}'(\alpha) dim_{\alpha}(\lambda) Jack_{\alpha}(\Lambda)}{dim_{\alpha}(\Lambda) Jack_{\alpha}(\lambda)}= \frac{c_{\lambda}(\alpha)}{c_{\Lambda}(\alpha)} \psi_{\Lambda/\lambda}'(\alpha).\] If $\lambda$ is chosen from Jack$_{\alpha}$ measure on partitions of size $n$, then after this transition $\Lambda$ is chosen from Jack$_{\alpha}$ measure on partitions of size $n+1$. Moreover there is a way of transitioning from a partition $\lambda$ of size $n$ to a partition $\tau$ of size $n-1$. The transition to $\tau$ occurs with probability \[ \frac{dim_{\alpha}(\tau) \psi_{\lambda/\tau}'(\alpha)}{dim_{\alpha}(\lambda)} = \frac{\psi_{\lambda/\tau}'(\alpha) c_{\lambda}'(\alpha)}{\alpha n c_{\tau}'(\alpha)}.\] If $\lambda$ is chosen from Jack$_{\alpha}$ measure on partitions of size $n$, then after this transition $\tau$ is chosen from Jack$_{\alpha}$ measure on partitions of size $n-1$.

	Now we can define analogs of the random variables $Y_1,\cdots,Y_n$ of Section \ref{charmart}. Let $\mu$ be a partition of $n$, and let $\lambda(1),\cdots,\lambda(n)$ be a sequence of random partitions generated by Kerov's growth process starting from $\lambda(1)$ equal to the only partition of size 1. Now define $Y_j^{(\alpha)} = \theta^{\lambda(j)}_{\mu-1^{n-j}}(\alpha)$, where this is understood to be 0 if $\mu$ has fewer than $n-j$ parts of size $1$, or if $j=0$.

\begin{theorem} \label{ismart2} For any partition $\mu$ of $n$, the sequence $(Y_1^{(\alpha)},\cdots,Y_n^{(\alpha)})$ is a martingale. \end{theorem}

\begin{proof} It is necessary to show that for all $\lambda$ of size $j$, \[ \sum_{|\Lambda|=j+1} \psi_{\Lambda/\lambda}'(\alpha) \frac{c_{\lambda}(\alpha)}{c_{\Lambda}(\alpha)} \theta^{\Lambda}_{\mu-1^{n-j-1}}(\alpha) = \theta^{\lambda}_{\mu-1^{n-j}}(\alpha).\] But this is clear from part 1 of Lemma \ref{crux}. \end{proof}

	Proposition \ref{condit2} is interesting because it shows that the martingale of Theorem \ref{ismart} has special properties. However it will not be needed in the sequel.

\begin{prop} \label{condit2} Let $\mu$ be a partition of $n$. Let $x_{\downarrow i}$ denote $(x)(x-1) \cdots (x-i+1)$ or $1$ in the case that $i=0$. Then for all $1 \leq j \leq n$, \[ E^{Y_n^{(\alpha)}}(Y_j^{(\alpha)}) = Y_n^{(\alpha)} \left( \frac{m_1(\mu)_{\downarrow n-j}}{n_{\downarrow n-j}}  \right).\] \end{prop} 

\begin{proof} In fact we show that \[ E^{\lambda(n)} (Y_j^{(\alpha)}) = Y_n^{(\alpha)} \left( \frac{m_1(\mu)_{\downarrow n-j}}{n_{\downarrow n-j}}  \right)\] which is stronger. 

	As the case $j=n$ is clear, we first treat the case $j=n-1$. Observe that
\[ E^{\lambda(n)}(Y_{n-1}^{(\alpha)}) =
\sum_{|\tau|=n-1}
\frac{P(\lambda(n-1)=\tau,\lambda(n)=\lambda)}{P(\lambda (n)=\lambda)}
\theta^{\tau}_{\mu-1}(\alpha). \] From the description of Kerov's transition mechanism before the proof of Theorem \ref{ismart2}, it is clear that
$\frac{P(\lambda(n-1)=\tau,\lambda(n)=\lambda)}{P(\lambda
(n)=\lambda)}$ is simply the probability of transitioning down from
$\lambda$ to $\tau$, which is $\frac{\psi_{\lambda/\tau}'(\alpha)
c_{\lambda}'(\alpha)}{\alpha n c_{\tau}'(\alpha)}$. Substituting this
into the expression for $E^{\lambda(n)}(Y_{n-1}^{(\alpha)})$ shows that the latter
is equal to \[ \sum_{|\tau|=n-1}
\theta^{\tau}_{\mu-1}(\alpha) \frac{\psi_{\lambda/\tau}'(\alpha)
c_{\lambda}'(\alpha)}{\alpha n c_{\tau}'(\alpha)}.\] The case $j=n-1$
now follows from part 2 of Lemma \ref{crux}.

	Next we indicate how to treat the case $j=n-2$, from which the general argument is completely clear. By considerations similar to those in the previous paragraph, $E^{\lambda(n)}(Y_{n-2}^{(\alpha)})$ is equal to \begin{eqnarray*} & &  \sum_{|\tau|=n-2}
\frac{P(\rho(n-2)=\tau,\rho(n)=\lambda)}{P(\rho(n)=\lambda)}
\theta^{\tau}_{\mu-1^2}(\alpha)\\
 & = & \sum_{|\tau|= n-2} \sum_{|\eta|=n-1}
\frac{P(\rho(n-2)=\tau,\rho(n-1)=\eta,\rho(n)=\lambda)}{P(\rho(n-1)=\eta,\rho(n)=\lambda)}\\
& & \cdot
\frac{P(\rho(n-1)=\eta,\rho(n)=\lambda)}{P(\rho(n)=\lambda)}
\theta^{\tau}_{\mu-1^2}(\alpha). \end{eqnarray*} Arguing as in the previous paragraph, one sees that \begin{eqnarray*}
& & \sum_{|\tau|=n-2} \frac{P(\rho(n-2)=\tau,\rho(n-1)=\eta,\rho(n)=\lambda)}{P(\rho(n-1)=\eta,\rho(n)=\lambda)} \theta^{\tau}_{\mu-1^2}(\alpha)\\
& = & \sum_{|\tau|=n-2} \frac{\psi_{\eta/\tau}'(\alpha) c_{\eta}'(\alpha)}{\alpha (n-1) c_{\tau}'(\alpha)} \theta^{\tau}_{\mu-1^2}(\alpha)\\
& = &  \left( \frac{m_1(\mu)-1}{n-1} \right) \theta^{\eta}_{\mu-1}(\alpha). \end{eqnarray*} The result now follows by the fact in the previous paragraph that \[ \sum_{|\eta|=n-1} \frac{P(\rho(n-1)=\eta,\rho(n)=\lambda)}{P(\rho(n)=\lambda)} \theta^{\eta}_{\mu-1}(\alpha) = \frac{m_1(\mu)}{n} \theta^{\lambda}_{\mu}(\alpha).\] \end{proof}

\begin{prop} \label{exactfor2} Let $\mu$ be any partition of $n$. Then \[ E(Y_j^{(\alpha)}-Y_{j-1}^{(\alpha)})^2=\left\{ \begin{array}{ll}
0 & \mbox{if \ $m_1(\mu)<n-j$}\\
\frac{\alpha^{n-l(\mu)} (j-1)!}{z_{\mu-1^{n-j}}} \left( n-m_1(\mu) \right) & \mbox{otherwise }                                                \end{array}   \right. \] \end{prop}

\begin{proof} If $m_1(\mu)<n-j$ the result is clear since then by definition
 $Y_j^{(\alpha)}=Y_{j-1}^{(\alpha)}=0$. As in the proof of Proposition
 \ref{exactfor}, the martingale property of
 $(Y_1^{(\alpha)},\cdots,Y_n^{(\alpha)})$ implies that \[
 E(Y_j^{(\alpha)}-Y_{j-1}^{(\alpha)})^2=
 E(Y^{(\alpha)}_j)^2-E(Y^{(\alpha)}_{j-1})^2.\] Now one uses part 1 of
 Lemma \ref{useful} to compute that \[ E(Y_j^{(\alpha)})^2=
 \sum_{|\tau|=j} \frac{\alpha^j j!}{c_{\tau}(\alpha)
 c_{\tau}'(\alpha)} (\theta^{\tau}_{\mu-1^{n-j}}(\alpha))^2 =
 \frac{\alpha^{n-l(\mu)} j!}{z_{\mu-1^{n-j}}}. \] Similarly one has that
 $E(Y_{j-1}^{(\alpha)})^2=0$ if $m_1(\mu) =n-j$ and computes that
 otherwise is equal to \[ \frac{\alpha^{n-l(\mu)}
 {(j-1)}!}{z_{\mu-1^{n-j+1}}} = \frac{m_1(\mu)-(n-j)}{j}
 \frac{\alpha^{n
  -l(\mu)} j!}{z_{\mu-1^{n-j}}}. \] \end{proof}

\section{``Transpositions''} \label{trans2}

	This section examines and sharpens the results of Section \ref{jack} in the case when $\mu=(2,1^{n-2})$. We call this partition the class of ``transpositions'' even though the algebraic structure of conjugacy classes is not present. As in Section \ref{trans}, there are certain phenomena which occur only for the partition $\mu=(2,1^{n-2})$. 

To begin note that there is an explicit formula (\cite{Mac},page 384)
	\[ \theta^{\lambda}_{(2,1^{n-2})}(\alpha) = \sum_i \left(
	\alpha {\lambda_i \choose 2} - {\lambda_i' \choose 2} \right),
	\] where $\lambda_i$ is the length of row $i$ of $\lambda$ and
	$\lambda_i'$ is the length of column $i$ of $\lambda$. From
	this formula it follows that if $\lambda$ is obtained by
	adding a box $x$ to $\mu$, then \[
	\theta^{\lambda}_{(2,1^{n-2})}(\alpha) -
	\theta^{\mu}_{(2,1^{n-3})}(\alpha) = c_{\alpha}(x)\] where
	$c_{\alpha}(x)$ is the $\alpha$-content of $x$ defined as
	$\alpha$(column number of x-1)-(row number of x-1).

	Consequently Theorem \ref{ismart2} implies that for any partition $\lambda(j)$ of $j$, the expected value of the $\alpha$-content of the box added to $\lambda$ in moving from $\lambda(j)$ to $\lambda(j+1)$ is 0. This combinatorial fact was first established by Kerov \cite{K4}, whose proof method was completely different. In the same paper Kerov also proved that for {\it any} $\lambda$ of size $j$, the variance of the $\alpha$-content of the box added to $\lambda$ is equal to $\alpha |\lambda|$. This is stronger than Proposition \ref{exactfor2}, which implies this for a random Jack$_{\alpha}$ distributed $\lambda$--but unlike Kerov's result has the merit of being true for arbitrary partitions $\mu$. 

	Theorem \ref{Lmoment2} is a recent result of Lassalle. 

\begin{theorem}
 \label{Lmoment2} (\cite{La}, Theorem 8.1) Let
 $s_{r,\alpha}(\lambda)$ be the $r$th moment of the $\alpha$-content of the box
 added to $\lambda(j)=\lambda$ when transitioning to
 $\lambda(j+1)$. Let \[ d_{\rho}(\lambda) = \prod_{i \geq 1} (\sum_{x \in
 \lambda} c_{\alpha}(x)^i)^{m_i(\rho)}. \] Then $s_{r,\alpha}(\lambda)$ is equal to \[ \alpha^{r} \sum_{i=0}^{ \lfloor r/2 \rfloor}
 \sum_{h=0}^{r-2i} \sum_{k=0}^{min(i,h)} \frac{1}{\alpha^{i+h}}
 (1-\frac{1}{\alpha})^{r-2i-h} {|\lambda|+i-1 \choose i-k}
 \sum_{|\rho|=h \atop l(\rho) \leq k} u_{ihk}^{\rho}(r)
 \frac{d_{\rho}(\lambda)}{z_{\rho}}, \] where the coefficients $
 u_{ihk}^{\rho}(r)$ are certain positive integers discussed in
 \cite{La}. \end{theorem}

	For instance one has that
\[ \begin{array}{ll}
s_{1,\alpha}(\lambda)=0 & \\
s_{2,\alpha}(\lambda)= \alpha |\lambda| & \\
s_{3,\alpha}(\lambda)= 2 \alpha d_{1,\alpha}(\lambda) +\alpha (\alpha-1) |\lambda|& \\
s_{4,\alpha}(\lambda)= 3 \alpha d_{2,\alpha}(\lambda) + 3 \alpha (\alpha-1) d_{1,\alpha}(\lambda) + \alpha^2 {|\lambda|+1 \choose 2} + \alpha(\alpha-1)^2 |\lambda|&
\end{array}\]

	In order to use Theorem \ref{Lmoment2} to analyze the Jack$_{\alpha}$ average of $s_{r,\alpha}(\lambda)$, it is necessary to upper bound the Jack$_{\alpha}$ average of $d_{\rho}(\lambda)$. This is done in Lemma \ref{rhobound}. 

\begin{lemma} \label{rhobound} Suppose that $\alpha \geq 1$. Then if $\rho$ is a partition of size $h$ with at most $k$ parts, there is a constant $C_{k,h,\alpha }$ (depending only on $k,h,\alpha$) such that if $\lambda$ is chosen from $Jack_{\alpha}$ measure on partitions of size $n$, then \[ E(d_{\rho}(\lambda)) \leq C_{k,h,\alpha} n^{k+\frac{h}{2}}.\] \end{lemma}

\begin{proof}
 From Lemma 6.6 of \cite{F2} and the fact that $\alpha$ is fixed, it
 follows that the Jack$_{\alpha}$ probability that $\lambda_1 \geq 2e
 \sqrt{\frac{n}{\alpha}}$ or $\lambda_1' \geq 2e \sqrt{\alpha n}$
 decays exponentially as a function of $\sqrt{n}$. Thus excluding this
 unlikely event, one sees that the $\alpha$-content of any of the $n$
 boxes of $\lambda$ is at most $2e \sqrt{\alpha n}$ in absolute
 value. The lemma follows easily. \end{proof}

	Next we obtain an upper bound on the moments
	$E(Y^{(\alpha)}_{j+1}-Y^{(\alpha)}_j)^{2r}$.

\begin{prop} \label{lastbound} Suppose that $\alpha \geq 1$. Then there is a constant $C_{2r,\alpha}$ (depending only on $2r$ and $\alpha$) such that $E(Y^{(\alpha)}_{j+1}-Y^{(\alpha)}_j)^{2r} \leq C_{2r,\alpha} j^r$ for all $j \geq 1$. \end{prop}

\begin{proof}
 By Theorem \ref{Lmoment2} and Lemma \ref{rhobound}, it follows that
 for fixed values of $i,h,k$, the Jack$_{\alpha}$ average of the
 corresponding term in Lassalle's formula for $s_{2r,\alpha}(\lambda)$ where $\lambda$ has size $j$
 is at most \[ C'_{2r,\alpha} j^{i-k} j^{k+\frac{h}{2}} \leq C'_{2r,\alpha} j^r. \] Here $C'_{2r,\alpha}$ is a constant depending only on
 $2r,\alpha$. The result follows. \end{proof}

	To close this section, we give a generalization of Proposition
	\ref{squaremart}.

\begin{prop} Let $(Y_1^{(\alpha)},\cdots,Y_n^{(\alpha)})$ be the martingale of Theorem \ref{ismart2} corresponding to the choice 
$\mu=(2,1^{n-2})$. Then \[ \left ((Y_1^{(\alpha)})^2,(Y_2^{(\alpha)})^2-\alpha {2 \choose 2},\cdots,(Y_n^{(\alpha)})^2-\alpha {n \choose 2} \right) \] is a martingale. \end{prop}

\begin{proof} The proof method is the same as in Proposition \ref{squaremart}; one replaces all occurrences of $c(x)$ by $c_{\alpha}(x)$ and ${j \choose 2}$ by $\alpha {j \choose 2}$. One also needs the facts that the expected value of $c(y)$ given $\mu$ is 0 and the expected value of $c(y)^2$ given $\mu$ is $\alpha |\mu|$. \end{proof} 

\section{Central limit theory for martingales} \label{CLT}

	This section uses martingale theory to prove a central limit theorem for $\frac{\theta^{\lambda}_{(2,1^{n-2})}(\alpha)}{\sqrt{\alpha {n \choose 2}}}$ under Jack$_{\alpha}$ measure. Namely we prove Theorem \ref{sharptrans2} (and hence also Theorem \ref{sharptrans} which is a special case).

	We apply the following recent result of Haeusler \cite{Ha}. For general background on central limit theorems for martingales, see the references in the introduction of \cite{Ha}, most notably \cite{Bol} and \cite{HH}.

\begin{theorem} \label{hau} (\cite{Ha}) Let the real valued random variables $X_{1},\cdots,X_{n}$ be a square integrable martingale difference sequence with respect to the $\sigma$-fields ${\it F}_{0} \subset {\it F}_{1} \cdots \subset {\it F_{n}}$. In other words we suppose that $E(X_j^2) < \infty$ and $E(X_j|{\it F}_{j-1})=0$ for all $j$. For $\delta>0$ let
\[ L_{n,2 \delta} = \sum_{j=1}^n E(|X_j|^{2+2 \delta}) \] and \[ N_{n,2 \delta} = E \left( \left| \sum_{j=1}^n E(X_j^2|{\it F}_{j-1})-1 \right|^{1+\delta} \right).\] Then there is a constant $C_{\delta}$ depending only on $\delta$ such that for all real $x_0$, \[ \left| P(X_1+\cdots+X_n \leq x_0) - \frac{1}{\sqrt{2 \pi}} \int_{-\infty}^{x_0} e^{-\frac{x^2}{2}} dx \right| \leq C_{\delta} (L_{n,2 \delta}+ N_{n,2 \delta})^{1/(3+2 \delta)}.\] \end{theorem}

	To prove Theorem \ref{sharptrans2} we will apply Theorem \ref{hau} with $X_j=\frac{Y_j^{(\alpha)}-Y_{j-1}^{(\alpha)}}{\sqrt{\alpha {n \choose 2}}}$ and $\mu=(2,1^{n-2})$. In this case (the Jack analog of the character ratio on transpositions), two major simplifications will occur: the quantity $N_{n,2 \delta}$ is in fact 0, and the quantity $L_{n,2 \delta}$ can be bounded using Proposition \ref{lastbound}. Bounding these quantities for other partitions seems to be a difficult problem, but as mentioned in the introduction the class of transpositions seems to have unique importance.

\begin{proof}
 (Of Theorem \ref{sharptrans2}) We apply Theorem \ref{hau} with
 $X_j=\frac{Y_j^{(\alpha)}-Y_{j-1}^{(\alpha)}}{\sqrt{\alpha {n \choose
 2}}}$ and $\mu=(2,1^{n-2})$. Let $\delta$ be a positive integer. By
 work of Kerov \cite{K4}, $s_{2,\alpha}(\lambda)=\alpha|\lambda|$ is
 independent of $\lambda$, so $N_{n,2 \delta}=0$. By Proposition
 \ref{lastbound}, \[ L_{n,2 \delta} \leq \left( \frac{1}{\sqrt{\alpha {n
 \choose 2}}} \right)^{2+2\delta} \sum_{j=1}^{n} C_{2 \delta+2,\alpha}
 j^{\delta+1} \leq \frac{C'_{2 \delta+2,\alpha}}{n^{\delta}}.\] Here
 $C'_{2 \delta+2,\alpha}$ is a constant depending only on $\delta$ and
 $\alpha$. Now observe that for any $s<\frac{1}{2}$, one can find
 $\delta$ large enough so that $\left(\frac{1}{n^{\delta}}
 \right)^{1/(3+ 2 \delta)} < n^{-s}$. \end{proof}

{\bf Remark:} The reader might be inclined to conjecture that any
martingale sequence (for instance ours) with $N_{n,2 \delta}=0$ and
the $X_j$ having uniformly bounded third moment would yield a central
limit theorem with convergence rate of order $n^{-1/2}$. Grams
\cite{Gr} shows that such sequences give central limit theorems with
convergence rate of order $n^{-1/4}$, but Bolthausen \cite{Bol} gives
examples showing that even with these hypotheses one can not in general
beat the $n^{-1/4}$ rate.

\section{Probabilistic proof of Burnside-Brauer} \label{probproof}

	In order to study the distribution of character ratios and their Jack analogs by Stein's method, the paper \cite{F2} showed how any representation of a finite group $G$ with real valued character can be used to construct a natural Markov chain on the set of irreducible representations of $G$. The purpose of this section is to show how this Markov chain gives a probabilistic proof of the following result of Burnside and Brauer (proved in \cite{Is}), in the special case that $\chi$ is real valued. 

\begin{theorem} \label{bb} (Burnside,Brauer) Let $\eta$ be the character of a faithful representation of a finite group $G$. Let $m$ be the number of distinct values assumed by $\eta$. Then every irreducible character of $G$ occurs as a constituent of some $\eta^j$, where $0 \leq j<m$. \end{theorem}

	In fact it will be shown that this result is closely related to a generalization of the following fact from algebraic graph theory.

\begin{theorem} \label{biggs} (\cite{Bi}) A connected graph with diameter $d$ has at least $d+1$ distinct eigenvalues. \end{theorem}

	To explain this, a dictionary is needed to go between representation theory and graph theory. This was developed in \cite{F2} and we recall it. Let $\eta$ be a real valued character of a finite group $G$, and let $<\phi,\theta>$ be the usual inner product on class functions of $G$ defined as $\frac{1}{|G|} \sum_{g \in G} \phi(g) \overline{\theta(g)}$. Then one can define a Markov chain $L_{\eta}$ on the set of irreducible representations of $G$ which transitions from $\lambda$ to $\rho$ with probability \[ L_{\eta}(\lambda,\rho) = \frac{dim(\rho)}{dim(\eta)dim(\lambda)} <\chi^{\rho},\eta \chi^{\lambda}>.\] Note that this transition probability is non-negative since $<\chi^{\rho},\eta \chi^{\lambda}>$ is the multiplicity of $\rho$ in the tensor product of $\lambda$ and the representation with character $\eta$. 

	Recall that a Markov chain on a finite set $X$ with transition probability $K(x,y)$ is said to be reversible with repect to a probability measure $\pi$ on $X$ if $\pi(x) K(x,y) = \pi(y) K(y,x)$ for all $x,y$. It is straightforward to see that this condition implies that $\pi$ is a stationary distribution
for $K$ (i.e. that $\pi(y) = \sum_x \pi(x)K(x,y)$ for all $y$). Moreover consider the space of real valued functions
$\ell^2(\pi)$ with the norm \[ ||f||_2 = \left( \sum_x |f(x)|^2
\pi(x) \right) ^{1/2}\] and let $K$ denote the Markov
operator on $\ell^2(\pi)$ defined by \[Kf(x) = \sum_y K(x,y)
f(y).\] Then the operator $K$ is self-adjoint so has real
eigenvalues $-1 \leq \beta_{min}=\beta_{|X|-1} \cdots \leq
\beta_0=1$. Moreover if $\{ \psi_i \}$ is an orthonormal basis of eigenfunctions with eigenvalues $\beta_i$, elementary linear algebra implies that that the chance that the chain started at $x$ is at $y$ after $r$ steps is $\sum_i \beta_i^r \psi_i(x) \psi_i(y) \pi(y)$.

\begin{theorem} \label{dict} (\cite{F2})
\begin{enumerate}
\item The transition probabilities of $L_{\eta}$ sum to 1, and the Markov chain $L_{\eta}$ is reversible with respect to the Plancherel measure $\pi$ of $G$.
\item The eigenvalues of the chain $L_{\eta}$ are indexed by conjugacy classes $C$ of $G$. The eigenvalue parameterized by $C$ is the character ratio $\frac{\eta(C)}{\eta(1)}$, and an orthonormal basis of eigenfunctions $\psi_C$ in $L^2(\pi)$ is defined by $\psi_C(\rho)= \frac{|C|^{1/2} \chi^{\rho}(C)}{dim(\rho)}$.
\end{enumerate}
\end{theorem}
 
	One can think of any reversible Markov chain on a finite set as inducing a graph structure on that set, with the unordered edge $(x,y)$ having weight $\pi(x) K(x,y)$. In the case at hand, the vertices of the graph are the irreducible representations of $G$, and an edge $(\lambda,\rho)$ is given weight $\pi(\lambda) L_{\eta}(\lambda,\rho)$, where $\pi$ is the Plancherel measure of $G$. 

	Now we give a proof of Theorem \ref{bb}, in the case that $\eta$ is real valued.

\begin{proof} Suppose that $\eta$ is real valued, so that Theorem \ref{dict} is applicable. By Theorem \ref{dict}, the number of distinct eigenvalues of $L_{\eta}$ is equal to the number of distinct character values of $\eta$. Moreover since $\eta$ is faithful, the eigenvalue 1 occurs with multiplicity one and it follows that the weighted graph associated to $L_{\eta}$ is connected (since otherwise each connected component could be used to construct an eigenvalue 1 eigenfunction supported on that component). Thus by a straightforward generalization of the proof of Theorem \ref{biggs} to graphs with edge weights, it follows that the diameter of the graph associated to $\eta$ is less than $m$, the number of distinct character values assumed by $\eta$.  

	To complete the proof it is sufficient to show that $\chi^{\rho}$ occurs as a constituent in $\eta^j$ if and only if there is a path of length $j$ from the trivial representation to $\rho$ in the weighted graph corresponding to $L_{\eta}$. But Theorem \ref{dict} gives a complete diagonalization of the Markov chain $L_{\eta}$, so by the sentence preceding the statement of Theorem \ref{dict}, it follows that the chance that $L_{\eta}$ transitions from the trivial representation to $\rho$ in $j$ steps is equal to \begin{eqnarray*}
& & \sum_{C} \left( \frac{\eta(C)}{\eta(1)} \right)^j |C| \frac{\chi^{\rho}(C)}{dim(\rho)} \frac{dim(\rho)^2}{|G|}\\
& = & \frac{dim(\rho)}{\eta(1)^j} \frac{1}{|G|} \sum_{g \in G} \eta^j(g) \chi^{\rho}(g). \end{eqnarray*} The result now follows since \[  \frac{1}{|G|} \sum_{g \in G} \eta^j(g) \chi^{\rho}(g) \] is the multiplicity of $\chi^{\rho}$ in $\eta^j$. \end{proof}

\section{Acknowledgements} The author was partially supported by National Security Agency grant MDA904-03-1-0049. The author thanks Michel Lassalle for correspondence.

\end{document}